\documentstyle[12pt]{article}
\newcommand{\const}{\mathop{\rm const}\limits}

\newcommand{\supp}{\mathop{\rm supp}\limits}

\begin{document}

\begin{center}

\vspace{3mm}

{\bf Lebesgue Spaces Norm Estimates for }\par

\vspace{4mm}

{\bf Fractional Integrals and Derivatives. } \\

\vspace{4mm}

 $ {\bf E.Ostrovsky^a, \ \ L.Sirota^b } $ \\

\vspace{4mm}

$ ^a $ Corresponding Author. Department of Mathematics and computer science, Bar-Ilan University, 84105, Ramat Gan, Israel.\\

E-mail: eugostrovsky@list.ru\\

\vspace{3mm}

$ ^b $  Department of Mathematics and computer science. Bar-Ilan University,
84105, Ramat Gan, Israel.\\

E-mail: sirota3@bezeqint.net \\

\vspace{4mm}
                    {\sc Abstract.}\\

 \end{center}

 \vspace{3mm}

  We study the problem estimation of classical Lebesgue-Riesz and Grand Lebesgue Norm for
 the fractional integrals and derivatives for the functions  from the classical Lebesgue-Riesz spaces
 as well as from the modified Besov's spaces. \par

  \vspace{4mm}

 \par

\vspace{3mm}


\vspace{4mm}

\section{Notations. Statement of problem.}

\vspace{4mm}

 "Fractional derivatives have been around for centuries  but recently they have
found new applications in physics, hydrology and finance", see  \cite{Meerschaert1}. \par
 Another applications: in the theory of Differential Equations are described in \cite{Miller1};
in statistics see in \cite{Bapna1}, \cite{Ostrovsky8}; see also \cite{Golubev1}, \cite{Enikeeva1}; in the theory of
integral equations etc. see in  the classical monograph \cite{Samko1}.  \par

\vspace{4mm}

  Let $ \alpha = \const \in (0,1); $ and let $ g = g(x), \  x \in R_+ $  be measurable numerical function. The fractional derivative
 of a Riemann - Liouville type of order $  \alpha: \  D^{\alpha}[g](x)  = g^{(\alpha)}(x) $  \cite{Riemann1}, \cite{Liouville1}
 is defined as follows: $ \Gamma(1-\alpha)  g^{(\alpha)}(x) =  $

$$
\Gamma(1-\alpha) \ D^{\alpha}[g](x)  =
\Gamma(1-\alpha) \ D^{\alpha}_x[g](x) \stackrel{def}{= } \frac{d}{dx} \int_0^x \frac{g(t) \ dt}{(x-t)^{\alpha}}. \eqno(1.1)
$$
see, e.g. the classical monograph of S.G.Samko, A.A.Kilbas and O.I.Marichev \cite{Samko1}, pp. 33-38; see also \cite{Miller1}.\par

 The case when $ \alpha \in (k, k + 1), \ k = 1,2,\ldots $ may be considered analogously, through the suitable derivatives of integer order. \par

 Hereafter $ \Gamma(\cdot) $ denotes the ordinary Gamma  function. \par

 We agree to take $ D^{\alpha}[g](x_0) = 0,  $ if at the point $  x_0 $ the expression $  D^{\alpha}[g](x_0)  $ does not exists. \par

Notice that the operator of the fractional derivative is non-local, if $ \alpha $ is not integer non-negative number.\par

 Recall also that the fractional  integral $  I^{(\alpha)}[\phi](x) = I^{\alpha}[\phi](x)  $ of a Riemann-Liouville type of
an order $ \alpha, 0 < \alpha < 1 $ is defined as follows:

$$
I^{(\alpha)}[\phi](x) \stackrel{def}{=} \frac{1}{\Gamma(\alpha)} \cdot \int_0^x \frac{\phi(t) \ dt}{(x-t)^{1 - \alpha}}, \ x,t > 0. \eqno(1.2)
$$
 It is known  (theorem of Abel, see \cite{Samko1}, chapter 2, section 2.1)
 that the operator $ I^{(\alpha)}[\cdot]  $ is inverse to the fractional derivative operator $ D^{(\alpha)}[\cdot],  $
at least  in the class of absolutely continuous functions. \par

 Another approach to the introducing of the fractional derivative, more exactly, the fractional Laplace operator
leads us to the using of Fourier transform

$$
 F[f](t) = \int_{R^d} e^{i (t,x) } \ f(x) \ dx
$$
in the space $  R^d, \ d = 1,2,\ldots:  $

$$
R_{\alpha,F}[f] := C_1(\alpha,d) F^{-1} \left[ |x|^{\alpha} F[f](x)  \right], \ 0 < \alpha < d,
$$
which leads us in turn up to multiplicative constant to the well-known Riesz potential

$$
R_{\alpha}[f] \stackrel{def}{=} \int_{R^d} \frac{f(y) \ dy}{ |x-y|^{d - \alpha}}, \ 0 < \alpha < d. \eqno(1.3)
$$

 Hereafter $ (t,x) = \sum_{m=1}^d t_m x_m, \ |x| = \sqrt{ (x,x) }, \ t,x \in R^d.  $ \par

\vspace{4mm}

 {\bf  We consider in this short article  the problem of Grand Lebesgue Norm estimation for
 the fractional integrals and derivatives for the functions  from the classical Lebesgue-Riesz spaces as well as from the Besov spaces.} \par

\vspace{3mm}

 Recall that the classical Lebesgue-Riesz  $  L(p) $ norm $  |f|_p  $ of a function $  f $ is defined by a formula

$$
|f|_p = \left[ \int_{R^d} |f(x)|^p \ dx  \right]^{1/p}, \ 1 \le p < \infty
$$
or correspondingly

$$
|f|_p = \left[ \int_{R_+} |f(x)|^p \ dx  \right]^{1/p}, \ 1 \le p < \infty.
$$

\vspace{4mm}

\section{Fractional integral estimate in the classical Lebesgue-Riesz norm.}

\vspace{4mm}

 Let for beginning $  x \in R^d, \ d = 1,2,\ldots, \ 0 < \alpha < d, \ p_+ = d/\alpha. $ Define for the value $ p = \const \in (1, p_+) $
the variable $ q $ as follows

$$
\frac{1}{q} \stackrel{def}{=} \frac{1}{p} - \frac{\alpha}{d}, \eqno(2.0)
$$
then $ q \in (d/(d - \alpha), \infty).  $ The relation (2.0) defines the value $  q  $ as a unique defined function on $  p: \ q = q(p) $
and conversely $ p = p(q). $\par

 Let us investigate in this section the inequality of a form

$$
|R_{\alpha}[f]|_q \le K_{R,\alpha}(p) \cdot |f|_p, \ f \in L_p(R^d) \eqno(2.1)
$$
or as a particular case

$$
| \Gamma(\alpha) \ I^{\alpha}[f]|_q \le K_{I,\alpha}(p) \cdot |f|_p, \ f \in L_p(R_+). \eqno(2.1a)
$$

The equality (2.0) is necessary and sufficient for the existence and finiteness of the "constant"  $ K_{R,\alpha}(p), $
see \cite{Ostrovsky10} - \cite{Ostrovsky13}, \cite{Stein1}, as well as the restriction $ 1 < p < d/\alpha. $
Obviously, in the case of the fractional integration (2.1a) $  d = 1 $ and hence $ 1/q = 1/p - \alpha, \ 0 < \alpha < 1, \ 1 < p < 1/\alpha. $ \par

\vspace{3mm}

 This classical problem goes back to Hardy and Littlewood, see \cite{Hardy0}, \cite{Hardy1},  \cite{Hardy2},   \cite{Hardy3}; more modern works
\cite{Adams1}, chapter 3, \cite{Adams2}, \cite{Adams3},  \cite{Adams4},
\cite{Harboure1}, \cite{Ostrovsky10} - \cite{Ostrovsky13}, \cite{Samko1}, p. 64-76 etc.\par

\vspace{3mm}

 We will understand  in the sequel in the capacity of the value $ K_{R,\alpha}(p)  $ its minimal value, namely

$$
K_{R,\alpha}(p) \stackrel{def}{=}
\sup_{ 0 \ne f \in L(p) } \left[\frac{|R_{\alpha}[f]|_q}{|f|_p} \right], \ 1/q = 1/p - \alpha/d, \ 1 < p < d/\alpha, \eqno(2.2)
$$
and correspondingly

$$
K_{I,\alpha}(p) \stackrel{def}{=}
\sup_{ 0 \ne f \in L(p) } \left[\frac{|I^{\alpha}[f]|_q}{|f|_p} \right], \ 1/q = 1/p - \alpha, \ 1 < p < 1/\alpha. \eqno(2.2a)
$$

 In order to formulate and prove a main result of this section, we need to introduce some notations. Put

$$
\Omega(d) = \frac{\pi^{d/2}}{\Gamma(d/2 + 1)} \ -
$$
be a volume of Euclidean unit $ d \ - $  ball. \par
 Further, the so-called maximal operator $ M[f](x), \ x \in R^d $ is defined as follows:

$$
M[f](x) \stackrel{def}{=} \sup_{r > 0} \left\{ \left( \Omega(d) \ r^d  \right)^{-1} \int_{B(x,r)} |f(x)| \ dx \right\},
$$
where $ B(x,r) $ is closed Euclidean ball in the whole space $  R^d  $ with center at the point $ x  $ and with radii  $ r,  \ r > 0. $ \par
 It is  known that

$$
|M[f]|_p \le S(d) \cdot \frac{p}{p-1} \cdot |f|_p, \ p \ge 1,
$$
where the finite "constant" $  S(d) $ is named as Stein's constant.\par

 The first upper estimation for the value $ S(d) $ was obtained in the classical book of
E.M.Stein  \cite{Stein1}, p. 173-188: $ S(d) \le 2 \cdot 5^d.$
It is proved in the article \cite{Hunt1}  that $ S(2) \le 2. $ In the next works of E.M.Stein \cite{Stein1} - \cite{Stein4}
was obtained consequently the following estimations for $ S(d):  \ S(d) \le c_1 \sqrt{d}, \ S(d) \le c_2 $
with some absolute constants $ c_1 $ and $ c_2. $ \par

\vspace{3mm}

{\bf Theorem 2.1.} We state under formulated above restrictions

$$
 K_{R,\alpha}(p) \le \frac{V_2(\alpha,d,p) \ \alpha^{-1} }{[(p-1) (1 - \alpha p) ]^{1 - \alpha/d   }}, \eqno(2.3)
$$
where

$$
 V_2(\alpha,d,p) = \Omega^{-1- \alpha/d }(d) \ p^{1 - 2 \alpha p/d} \ d^{ 1 + ( 1 - \alpha p  )/d   }  \
  S^{1 - \alpha p }(d) \in (0,\infty)  \eqno(2.3a)
$$
is continuous function relative the variable $  p  $ on the {\it closed} segment $  p \in [1, \ 1/\alpha];   $ \par

\vspace{3mm}

$$
 K_{I,\alpha}(p) \ge \frac{V_1(\alpha,p)}{[(p-1) (1 - \alpha p)]^{1 - \alpha }}, \eqno(2.4)
$$
where the "constant"  function  $  p \to V_1(\alpha,p) $ is also strictly positive continuous on the {\it closed} segment
$ p \in [1, 1/\alpha]. $\par

\vspace{3mm}

{\bf Proof of the upper bound.}\par

 It is sufficient to follow the book of  D.R.Adams \ \cite{Adams1}, chapter 3 and make accuracy computations.\par
 Some details. We write

$$
R_{\alpha}[f](x) = J_1 + J_2, \  J_1 = J_1(x) =\int_{ |x-y| < \delta } \frac{f(y) \ dy}{|x-y|^{d-\alpha} },
$$

$$
 J_2 = J_2(x) = \int_{ |x-y| \ge \delta } \frac{f(y) \ dy}{|x-y|^{d-\alpha} }, \ \delta > 0.
$$

 We use the H\"older's inequality for the $ J_2(\cdot) $ estimation:

$$
J_2(x) \le \Omega^{1 - 1/p}(d) \cdot \left( \frac{d-\alpha p}{p-1}  \right)^{-1 + 1/p} \cdot \delta^{\alpha - d/p} \cdot |f|_p.
$$

  Further, we apply  lemma 3.1.4 from the book \cite{Adams1}:

$$
\int_{|x-y|< \delta} \frac{\mu(dy)}{|x-y|^{d - \alpha}} = (d-\alpha) \int_0^{\delta} \frac{\mu(B(x,r)) \ dr}{r^{ 1 + d - \alpha  }} +
\frac{\mu(B(x,\delta))}{\delta^{d - \alpha}}. \eqno(2.5)
$$
 Here $  \mu(\cdot) $  is arbitrary Borelian measure in the space $ R^d;  $ indeed, in the considered case $ \mu(dy) = |f(y)| dy.  $ \par
Then

$$
\mu(B(x,\delta))= \int_{B(x,\delta)} |f(y) | \ dy =  \Omega(d) \ \delta^d \cdot \frac{1}{ \Omega(d) \ \delta^d } \int_{B(x,\delta)} |f(y)| \ dy \le
$$

$$
\Omega(d) \ \delta^d  \ M[f](x).
$$
 It remains to substitute into (2.5), apply the elementary equality
$$
\min_{\delta > 0} \left(  A \delta^{\alpha} + B \delta^{\alpha - d/p}  \right) =
\frac{d}{d - \alpha p} \cdot \left( \frac{d-\alpha p}{\alpha p} \right)^{\alpha p/d} \cdot A^{1 - \alpha p/d} \cdot B^{\alpha p /d}.
$$
and use further the integration over $ x, \ x \in R^d. $ \par

\vspace{4mm}

{\bf Proof of the lower bound.}\par

\vspace{3mm}

{\bf A. Case $ p \to 1+0. $  }\\

\vspace{3mm}

 Let us choose the following test function

 $$
 f_0(x) = x^{-1} \ I(x > 1).
 $$
Hereafter $  I(x \in A) = 1,  $ if $  x \in A $ and  $  I(x \in A) = 0, \ x \notin A.  $ \par

 We derive consequently

$$
|f_0|^p_p = \int_1^{\infty} x^{-p} \ dx = (p-1)^{-1}, \ p > 1;
$$

$$
|f_0|_p = (p-1)^{-1/p}, \ p > 1.
$$

 Further, when $ x > 1 $ and $ x \to \infty, \ q \to 1/(1 - \alpha)+ $

$$
g_0(x):= \int_0^x \frac{f_0(y) \ dy}{|x-y|^{1 - \alpha}} dy = \int_1^x \frac{y^{-1} \ dy}{|x-y|^{1 - \alpha}} dy =
$$

$$
x^{\alpha - 1} \  \int_{1/x}^1 \frac{z^{-1} \ dz}{|1 - z|^{ 1 - \alpha}} \sim x^{\alpha - 1} \ \ln x;
$$

$$
|g_0|_q^q \sim \int_1^{\infty} x^{(\alpha - 1)q} \ \ln^q (x) \ dx =
\frac{\Gamma(q + 1)}{ \left[ q(1-\alpha) - 1 \right]^{q + 1}};
$$

$$
|g_0|_q \asymp \left(q - 1/(1 - \alpha) \right)^{-1 - 1/q} \asymp (p-1)^{-1 - 1/q}.
$$
Therefore

$$
\frac{|g_0|_q}{|f_0|_p} \asymp \frac{(p-1)^{-1 - 1/q}}{(p-1)^{-1/p}} = (p-1)^{-1 - 1/q + 1/p } =
(p-1)^{-1 +\alpha}. \eqno(2.6)
$$

\vspace{3mm}

{\bf B. Case $ p \to 1/\alpha -0. $  }\\

\vspace{3mm}

 Let us choose now the following test function

$$
h_{\Delta} = h_{\Delta}(x) := x^{-\alpha} \ |\ln x|^{\Delta} \ I(0 < x < 1/e), \ \Delta = \const > 0. \eqno(2.7)
$$
 Further we will take $  \Delta \to 0+. $\par

 We have for the values $ x \to 0+, \ p \in (1, 1/\alpha), \ p \to 1/\alpha - 0  $ and correspondingly $ q \to \infty $

$$
|h_{\Delta}|_p^p = \int_0^{1/e} x^{-\alpha p} \ |\ln x|^{\Delta p} \ dx=
(1 - \alpha p)^{ -1 - \Delta p } \int_{1 - \alpha p}^{\infty}z^{\Delta p} \ e^{-z} \ dz \sim
$$

$$
(1 - \alpha p)^{ -1 - \Delta p }  \ \Gamma(\Delta p + 1); \hspace{6mm} |h_{\Delta}|_p \sim
\frac{\Gamma^{1/p}(\Delta p + 1)}{(1 - \alpha p)^{ \Delta + 1/p }}; \eqno(2.8)
$$

$$
r_{\Delta}(x) :=\Gamma(\alpha) I^{\alpha}[h_{\Delta}](x) = \int_0^{1/e} \frac{y^{-\alpha} |\ln y|^{\Delta}}{|x-y|^{1 - \alpha}} \ dy \sim
|\ln x|^{\Delta + 1};
$$

$$
|r_{\Delta}|_q \sim e^{-1} \ \alpha^{-\Delta - 1} \ (1 - \alpha p)^{- \Delta - 1}. \eqno(2.9)
$$

 Therefore

$$
\frac{|r_{\Delta}|_q}{|h_{\Delta}|_p} \sim \frac{e^{-1} \ \alpha^{-1}}{(1 - \alpha p)^{1 - \alpha}}. \eqno(2.10)
$$

\vspace{3mm}

{\bf C. General case $ (p-1)( 1/\alpha -p) \to 0. $  }\\

\vspace{3mm}

  Since the functions $ f_0(\cdot)  $ and $ r_{\Delta}(\cdot)  $ are disjoint:

$$
 f_0(x) \cdot r_{\Delta}(x)   = 0,
$$
we conclude

$$
|f_0(x) + r_{\Delta}(x) |_p^p =  |f_0(x)|_p^p + | r_{\Delta}(x) |_p^p.
$$

 Choosing ultimately in the capacity of the test function $ t(x)=f_0(x) + h_{\Delta}(x),  $ we get to
the second assertion of theorem 2.1.\par

\vspace{3mm}

 As a slight consequence: \\

{\bf Corollary 2.1.} If $ 1 < p < 1/\alpha, \ \alpha = \const \in (0,1) $ then

\vspace{3mm}

$$
 \frac{V_1(\alpha,p)}{[(p-1) (1 - \alpha p)]^{1 - \alpha }} \le K_{I,\alpha}(p) \le
 \frac{V_2(\alpha,1,p) \ \alpha^{-1} }{[(p-1) (1 - \alpha p)]^{1 - \alpha }}. \eqno(2.11)
$$

\vspace{4mm}

\section{Fractional integral estimate in the Grand Lebesgue norm.}

\vspace{4mm}

 Let  $ (X, Q, \nu )  $ be a measurable space with non - trivial sigma finite measure $  \nu, $ and let also $ \psi = \psi(q), \ s_1 \le q < s_2,
 \ 1 \le s_1 < s_2 \le \infty  $  be continuous on the open interval $   ( s_1, s_2 )  $ bounded from below  function.  By definition, a
Grand Lebesgue Space (GLS) $ G\psi = G\psi(s_1,s_2) $ over our triplet $ (X, Q, \nu )  $ consists  on all the measurable functions $  f: X \to R $  with
finite norm

$$
||f||G\psi \stackrel{def}{=} \sup_{q \in (s_1,s_2)} \left[  \frac{|f|_q}{\psi(q)}  \right]. \eqno(3.1)
$$
 Hereafter

 $$
 |f|_q  = \left[ \int_X |f(x)|^q \ \nu(dx)  \right]^{1/q}
 $$
and we will denote $ (s_1, s_2) = \supp \psi. $ \par
 The detail investigation  of these spaces see, e.g. in \cite{Fiorenza1},  \cite{Fiorenza2},
\cite{Iwaniec1}, \cite{Iwaniec2}, \cite{Kozachenko1},  \cite{Ostrovsky1}, \cite{Ostrovsky7}, \cite{Ostrovsky8} etc.\par

 The finiteness of Grand Lebesgue norm $  ||f||G\psi  < \infty $  in the case when $ s_2 = \infty $ implies in particular the exponential
decrease the tail $ T_f(u) = \nu \{ x: |f(x)| > u  \}, \ u \to \infty $ function of $  f.  $  \par

 The fundamental function $ \phi(G\psi, \delta), \delta > 0  $ for this spaces is defined as follows

$$
\phi(G\psi, \delta) = \sup_{p \in \supp \psi} \left[  \frac{\delta^{1/p}}{\psi(p)} \right].
$$

 This function play a very important role in the theory of operator's interpolation, theory of Fourier series etc. \cite{Bennet1}. \par

 In the considered here problems $ X = R^d  $ or $ X = R_+ $ with Lebesgue measure $ \nu; \ \nu(dx) = dx. $ \par

\vspace{3mm}

The set of all such a functions with support $ \supp(\psi) = (s_1, s_2) $ will be denoted by
$ G\Psi(s_1, s_2). $ \par

 These spaces are rearrangement invariant; and are used, for example, in the theory
of Probability, theory of Partial Differential Equations, Functional Analysis, theory of
Fourier series, Martingales, Mathematical Statistics, theory of Approximation etc.\par

 Notice that the classical Lebesgue-Riesz spaces $ L_p $ are extremal case of Grand Lebesgue
Spaces, see \cite{Ostrovsky16},   \cite{Ostrovsky17}. \par

 Let a function $ f : X \to R $ be such that

$$
\exists s_1, s_2, \ 1 \le s_1 < s_2 \le \infty: \ \forall p \in (s_1, s_2) \ \Rightarrow |f|_p < \infty.
$$

 Then the function $ \psi = \psi(p), \ s_1 < p < s_2  $ may be {\it naturally} defined by the following way:
$ \psi(p)  := |f|_p. $

\vspace{3mm}

Let now the (measurable) function $ f : X \to R $ be such that  $ f \in G\psi $ for some $ \psi(\cdot) $ with support
$ \supp \psi(\cdot) = (s_1,s_2)  $ for which $  1 \le s_1 < s_2 \le d/\alpha.$ Of course, the function $  \psi(p) $ may be
picked as a natural function for the function $  f(\cdot): $

$$
\psi(p) = \psi^{(f)}(p) := |f|_p,
$$
if it is finite for $  p \in (s_1, s_2). $\par
 We define a new $ \psi $ − function, say $ \psi_K = \psi_K(q) $ as follows.

$$
\psi_{K,R,\alpha}(q) = K_{R,\alpha}(p(q))\cdot \psi(p(q)), \ p \in (s_1,s_2).
$$

 Recall that the variable $ p $ and $  q  $ are closely related by the equality (2.2), which defined the variable
$  p  $ as unique function on $  q,  \ p = p(q); $  the conversely is also true.

\vspace{3mm}

{\bf Theorem 3.1.} We propose under formulated conditions, in particular under conditions of theorem 2.1

$$
||R_{\alpha}[f]||G\psi_{K,R,\alpha} \le 1 \cdot ||f||G\psi, \eqno(3.2)
$$
where the constant "1" in (3.2) is the best possible.\par

\vspace{3mm}

{\bf Proof. Upper bound.} \\

Let further in this section $   p \in (s_1, s_2), $  where $ 1 \le s_1 < s_2 \le d/\alpha.  $
We can and will suppose without loss of generality
$ ||f||G\psi = 1. $ Then $ |f|_p \le \psi(p),  p \in (s_1, s_2). $
We conclude after substituting into the inequality (2.3)

$$
|R_{\alpha}[f]|_q \le K_{R,\alpha}(p) \cdot \psi(p)= \psi_{K,R,\alpha}(q) = \psi_{K,R,\alpha}(q) ||f||G\psi.\eqno(3.3)
$$
The inequality (3.2) follows from (3.3) after substitution $ p = p(q). $ \par

\vspace{3mm}

{\bf Proof. Exactness.}\\
The exactness of the constant ”1” in the proposition (3.2) follows from the theorem 2.1
in the article \cite{Ostrovsky17}.\par

\vspace{4mm}

\section{Fractional derivative  of indicator function, with consequences.}

\vspace{4mm}

 The following example is computed and applied in \cite{Ostrovsky9}. Let here $ X = (0,1); $ define the function

$$
g_h(x) = I(h < x), \ x > 0, \ 0 < h = \const < 1,\eqno(4.0)
$$
 then

$$
g_h^{(\alpha)}(x) =  \frac{1}{\Gamma(1-\alpha)} \cdot I(h < x) \cdot (x - h)^{-\alpha}, \ \alpha = \const \in (0,1). \eqno(4.1)
$$

 It is verified also in \cite{Ostrovsky9} that

$$
I^{\alpha}  \left[ g_h^{(\alpha)} \right] (x)  =  I(h < x)  =   g_h(x). \eqno(4.2)
$$

 It is easily  to estimate $ |g_h|_p = h^{1/p},  $

$$
 \Gamma(1 - \alpha)\ |D^{\alpha} g_h(\cdot)|_q \asymp (1 - \alpha q)^{-1/q} \cdot h^{1/q - \alpha}, \ 1 \le p < \infty, 1 \le q < 1/\alpha.
$$

\vspace{3mm}

 We investigate in this section fractional derivative for a more general indicator function

$$
g_{h_1, h_2}(x) = I(h_1 < x < h_2), \ 0 < h_1 < h_2 < 1; \ \Delta := h_2 - h_1. \eqno(4.3)
$$
 Evidently,

$$
g_{h_1, h_2}(x) = g_{h_1}(x) - g_{h_2}(x),
$$
therefore

$$
\Gamma(1-\alpha) \ g_{h_1, h_2}^{(\alpha)}(x) = I(x > h_1) (x - h_1)^{-\alpha} - I(x > h_2) (x - h_2)^{-\alpha}.
$$

 We deduce after some computations for the values $ p: \ 1 \le p < 1/\alpha $

$$
\Delta^{1/p - \alpha} (1 - \alpha p)^{-1/p} \le |\Gamma(1-\alpha) \ g_{h_1, h_2}^{(\alpha)}(\cdot)|_p \le
3 \ \Delta^{1/p - \alpha} (1 - \alpha p)^{-1/p}. \eqno(4.4)
$$

  Let us estimate the fractional derivative $ D^{\alpha}[g_{h_1,h_2}] $ in the Grand Lebesgue Space norm.
Denote

$$
\psi_{\alpha}(p) =(1 - \alpha p)^{-1/p}, \ 1 < p < 1/\alpha.
$$
 Let also $  \zeta = \zeta(p), \ 1 < p < 1/\alpha  $ be any function from the set $ G\Psi(1, 1/\alpha),  $ then the product
function

$$
\theta(p) = \psi_{\alpha}(p) \cdot \zeta(p)
$$
belongs also the set $ G\Psi(1, 1/\alpha). $  \par

\vspace{3mm}

{\bf Proposition 4.1.}

$$
|| g^{(\alpha)}_{h_1,h_2} ||G\theta \le 3 \ \Delta^{- \alpha} \ \phi(G\zeta, \Delta), \ \Delta = h_2 - h_1. \eqno(4.5)
$$

\vspace{3mm}

{\bf Proof.} The right-hand side of the inequality (4.4) may be rewritten as follows.

$$
\Delta^{\alpha} \ | g^{(\alpha)}_{h_1,h_2}|_p \le 3 \ \psi_{\alpha}(p) \ \Delta^{1/p},
$$
or equally

$$
\Delta^{\alpha} \frac{ | g^{(\alpha)}_{h_1,h_2}|_p}{\theta(p)} \le 3 \  \frac{\Delta^{1/p}}{\zeta(p)}.
$$
 We deduce taking supremum from both the sides of the last inequality using the direct definition of Grand Lebesgue  norm and
fundamental function for these spaces

$$
 \Delta^{\alpha} \ ||g^{(\alpha)}_{h_1,h_2}||G\theta \le 3 \ \phi(G\zeta, \Delta),
$$
Q.E.D.\par

\vspace{3mm}

{\bf Definition 4.1.} A measurable function $  f: (0,b) \to R $  will named {\it very simple } with step $  h, 0 < h < 1,  $
write $ f \in VS(h), $ if it has a form

$$
f(x) = \sum_{k=1}^n c_k I_{A_k}(x), \eqno(4.6)
$$
where  $  A_k $ are pairwise disjoint segments  of a form $ A_k = (h_1(k), h_2(k)) $ with $  h_2(k) - h_1(k) = h = \const. $\par

 On the other words, $ f(\cdot) $ is spline of zero order with constant step, stepwise function.\par

 Obviously, if  $ f \in VS(h), $  then

$$
|f|_p = h^{1/p} \ \left[ \sum_k |c_k|^p  \right]^{1/p} = h^{1/p} \ |\vec{c}|_p, \ p \ge 1. \eqno(4.7)
$$

 But in the next pilcrow we impose  more strictly restriction  $ 1 \le p < 1/\alpha,  $ where again $  0 < \alpha < 1. $  We conclude
using triangle inequality and  inequality (4.4) \\ $ f \in VS(h) \ \Rightarrow $

$$
\Gamma(1 - \alpha) |f^{(\alpha)}|_p \le 3 \ h^{1/p - \alpha} \ (1-\alpha p)^{-1/p} \sum_k |c_k| =
$$

$$
 3 \ h^{1/p - \alpha} \ (1-\alpha p)^{-1/p}  |\vec{c}|_1 =   3 \ h^{1/p - \alpha -1} \ (1-\alpha p)^{-1/p} |f|_1. \eqno(4.8)
$$

\vspace{3mm}

 We obtain analogously and under conditions of proposition 4.1 \\

\vspace{3mm}

{\bf Proposition 4.2.} Let again  $ f \in VS(h), $ then

$$
|| f^{(\alpha)} ||G\theta \le 3 \ h^{- \alpha -1} \ \phi(G\zeta, h) \ |f|_1, \ h = h_2(k) - h_1(k) = \const. \eqno(4.9)
$$

 The last estimate may be used perhaps for numerical computation of fractional  derivatives via spline approximation.\par

 The condition $ h_2(k) - h_1(k) =  h =\const $ may be easily replaced to the following:

$$
 0 < C_1 < \frac{h_2(k)}{h_1(k)}  < C_2, \ C_1, C_2 = \const > 0.
$$

\vspace{4mm}

\section{Fractional derivative estimate.}

\vspace{4mm}

 The problem of norm estimation for fractional derivative is more complicated. \par

  Note first of all the definition (2.2a) common with estimation (2.4)  may be rewritten at least for absolutely
 continuous functions $  \{  f \} $ such that $ f(0) = 0 $ as follows:

$$
 |f|_q \le  K_{I,\alpha}(p) \cdot |D^{\alpha} \ f|_p, \ \alpha \in (0,1), \ 1/q = 1/p - \alpha, \ 1 < p < 1/\alpha \ -
$$
Sobolev's inequality for fractional derivatives. \par

  Let again $ \alpha = \const \in (0,1); $ and let $ f = f(x), \  x \in (0,b), \
 0 < b = \const \le \infty $  be measurable numerical function. The fractional derivative
 of a Riemann-Liouville type of order $  \alpha: \  D^{\alpha}[f](x)  = f^{(\alpha)}(x) $   is written, e.g. in (1.1). \par

 The next equality

$$
\Gamma(1 - \alpha) \ D^{\alpha}[f](x) = x^{-\alpha} f(x) + \alpha \int_0^x \frac{f(x) - f(t)}{(x - t)^{1 + \alpha}} \ dt
\stackrel{def}{=}
$$

$$
x^{-\alpha} f(x)  + U_{\alpha}[f](x), \eqno(5.0)
$$
which defines the so-called Marchaud fractional derivative, is proved e.g. in \cite{Samko1}, p. 220 - 229.\par

 Denote as ordinary by $ \omega(f,\delta)_p  $ the $  L_p  $ module of continuity of the function $  f:  $

$$
\omega(f,\delta)_p \stackrel{def}{=} \sup_{h: |h| \le \delta} |f(\cdot + h) - f(\cdot)|_p,
$$
where $  f(x) := 0, $  if $ \ x < 0 $ or if $ x > b  $ in the case when $ b < \infty. $  \par

 We introduce the following modification of the classical Besov's norm $ ||f||B^{(\alpha)}_p $
 and correspondent spaces $ \ B^{(\alpha)}_p $  as follows

$$
||f||B^{(\alpha)}_p \stackrel{def}{=} | x^{-\alpha} \ f(x) |_p + \alpha \int_0^b t^{-1 - \alpha} \omega(f,t)_p \ dt, \eqno(5.1)
$$

$$
f(\cdot) \in B^{(\alpha)}_p \ \Leftrightarrow ||f||B^{(\alpha)}_p  < \infty.
$$

\vspace{3mm}

{\bf Theorem 5.1.}

$$
\forall p \in (1,1/\alpha) \ \Rightarrow \sup_{0 \ne f \in B^{(\alpha)}_p  } \left\{ \frac{|D^{\alpha}[f]|_p}{||f||B^{(\alpha)}_p} \right\} =
\frac{1}{\Gamma(1 - \alpha)}. \eqno(5.2)
$$

\vspace{3mm}

{\bf Proof.} The inequality

$$
| \ \Gamma(1 - \alpha) \ D^{\alpha}[f] \ |_p  \le ||f||B^{(\alpha)}_p, \ f(\cdot) \in B^{(\alpha)}_p, \eqno(5.3)
$$
follows immediately from the representation (5.0), see \cite{Samko1}, p. 220-229.  In order to ground the lower bound for the
fraction

$$
r_{\alpha}(p):= \sup_{0 \ne f \in B^{(\alpha)}_p  } \left\{ \frac{|D^{\alpha}[f]|_p}{||f||B^{(\alpha)}_p} \right\}, \eqno(5.4)
$$
it is sufficient  to consider the following example (counterexample) with the value $  b = 1; $

$$
f_0(x) = g_{h_1, h_2}(x) = I(h_1 < x < h_2), \ h_1,h_2 = \const, \ 0 < h_1 < h_2 < 1.
$$
 Denote $ \Delta = h_2 - h_1,   $ then $ \Delta \in (0,1) $ and $ \Delta \to 1,  $  if $ h_1 \to 0+, \ h_2 \to 1-.  $
We estimate taking into account the restriction $ 1 \le p < 1/\alpha:  $

$$
||f_0||^p B^{(\alpha)}_p\ge \int_{h_1}^{h_2} x^{-\alpha p} \ dx = \frac{h_1^{1 - \alpha p}  - h_2^{1 - \alpha p}}{1 - \alpha p} \sim
\frac{\Delta^{1 - \alpha p}}{1 - \alpha p},
$$

$$
||f_0|| B^{(\alpha)}_p \ge  \frac{\Delta^{1/p - \alpha }}{(1 - \alpha p)^{1/p}}. \eqno(5.5)
$$

 The correspondent estimate for fractional derivative $ \Gamma(1 - \alpha) \ | f_0|_p $ is obtained in (4.4). Substituting into
expression (5.4), we  get to the proposition (5.2) of theorem 5.1.\par

\vspace{4mm}

{\bf  Remark 5.1.} Note that the inverse inequality for  (5.3), i.e. the inequality for arbitrary function $ f: R_+ \to R  $
of the form

$$
 ||f||B^{(\alpha)}_p  \le \tilde{K}_{I,\alpha}(p) \ \Gamma(1 - \alpha) \ | \ D^{\alpha}[f] \ |_p,  \ 0 < \alpha < 1, \
 \tilde{K}_{I,\alpha}(p) < \infty
$$
is not true for any number $  p, \ p \ge 1. $ Indeed, we put $ f_{\alpha}(x) = x^{\alpha - 1}; $ then
$   D^{\alpha}[f_{\alpha}] = 0, $ despite $  ||f_{\alpha}||B^{(\alpha)}_p > 0. $ \par

\vspace{3mm}

  Let us estimate the fractional derivative $ D^{\alpha}[f] $ in the Grand Lebesgue Space norm.
Suppose that there exists a value $ \beta, \ 1 \le \beta \le 1/\alpha  $ such that

$$
\psi^{(\beta)}(p) := |f||B^{(\alpha)}_p < \infty, \ 1 < p < \beta.
$$
 On the other words, the function $ p \to \psi^{(\beta)}(p)  $ is  Besov-Grand Lebesgue Spaces natural function for the function
$  f(\cdot). $ \par

\vspace{3mm}

{\bf Proposition 5.1. }

$$
|| \  D^{\alpha}[f] \ ||G\psi^{(\beta)} \le 1/\Gamma(1 - \alpha). \eqno(5.6)
$$

\vspace{3mm}

{\bf   Proof} is very elementary.  We use the equality (5.3) for the values $ p: 1 < p < \beta $

$$
| \ \Gamma(1 - \alpha) \ D^{\alpha}[f] \ |_p  \le ||f||B^{(\alpha)}_p  \le \psi^{(\beta)}(p),  \eqno(5.7)
$$
which is equivalent to the assertion of proposition (5.1.) \par

\vspace{4mm}

\section{Multidimensional case.}

\vspace{3mm}


 Let $  \alpha, \beta  = \const $ be two numbers such that $  0 < \alpha, \beta < 1; $
 The partial mixed fractional derivative $ D^{ \alpha, \beta}_{x,y}[G](x,y)   $
again of  Riemann-Liouville type
of order $ (\alpha, \beta) $ of a function $ G(\cdot, \cdot) $ at the positive points $ (x,y)  $ is defined as  follows:

$$
G^{(\alpha, \beta)}(x,y) =  D^{ \alpha, \beta}_{x,y}[G](x,y)  \stackrel{def}{=}   D^{\alpha}_x D^{\beta}_y [G] =
\frac{1}{\Gamma(1 - \alpha)} \frac{1}{\Gamma(1 - \beta)}\times
$$

$$
\frac{\partial^2}{ \partial x \partial y } \int_0^x \int_0^y \frac{G(t,s)\ dt \ ds}{ (x-t)^{\alpha} (y-s)^{\beta} }, \eqno(6.1)
$$
see, e.g. \cite{Samko1}, chapter 24. We put as before $ D^{ \alpha, \beta}_{x,y}[G](x,y) = 0 $ if at the point (in the plane) $ (x,y) $
the expression (4.2) for $ D^{ \alpha, \beta}_{x,y}[G](x,y) $ does not exists. \par

 Note that in general case $  D^{\alpha}_x D^{\beta}_y [H] \ne D^{\beta}_y D^{\alpha}_x [H],  $ but if the function $ G = G(x,y) $
is factorable: $ H(x,y) = g_1(x) g_2(y)  $ and both the functions  $ g_x(\cdot)  $ and $ g_2(\cdot) $ are "differentiable" at the
points $  x  $ and $  y  $ correspondingly:

$$
\exists D^{\alpha}[g_1](x), \hspace{6mm} \exists D^{\beta}[g_2](y),
$$
then really

$$
D^{\alpha}_x D^{\beta}_y [H] = D^{\beta}_y D^{\alpha}_x [H] = D^{\alpha}_x[g_1](x) \cdot D^{\beta}_y [g_2](y).
$$

 On the other words the (linear) operator $ D^{ \alpha, \beta}_{x,y}[\cdot] $ is the tensor product of the one - dimensional
fractional derivatives

$$
 D^{ \alpha, \beta}_{x,y}  =  D^{ \alpha}_{x}  \otimes   D^{\beta}_{y}.
$$

 Analogously be defined the partial mixed fractional integrals  $ I^{\alpha,\beta}[G], $ where $  0 < \alpha,\beta < 1: $

$$
 I^{ \alpha, \beta}_{x,y}  =  I^{ \alpha}_{x}  \otimes   I^{\beta}_{y}, \eqno(6.2)
$$
or in detail for the function $ G = G(x,y), \ x \in (0,b_1), \ y \in (0,b_2), \ 0 < b_{1,2} = \const \le \infty, $

$$
G^{(\alpha,\beta)}(x,y) =
I^{ \alpha, \beta}_{x,y} [G](x,y)  \stackrel{def}{=} \frac{1}{\Gamma(\alpha)} \
\frac{1}{\Gamma(\beta)} \int_0^x \int_0^y \frac{G(t,s)}{(x-t)^{1-\alpha} \ (y-s)^{1 - \beta}} dt \ ds. \eqno(6.3)
$$

 Recall that the so - called {\it  mixed Lebesgue - Riesz  } norm  $ ||f||_{p_1,p_2}, \ \\ 1 \le p_1, p_2 < \infty $ for a function  $  f(x,y) $
is defined by a formula

$$
 ||f||_{p_1,p_2} \stackrel{def}{=} \left\{ \int_0^{b_2}  \left[ \int_0^{b_1} |f(x,y)|^{p_1} \ dx   \right]^{p_2/p_1}  dy \right\}^{1/p_2}.
 \eqno(6.4)
$$

\vspace{3mm}

  Suppose $ 1 < p_1 < 1/\alpha, \ 1 < p_2 < 1/\beta  $ and define as before

$$
\frac{1}{q_1} = \frac{1}{p_1} - \alpha, \hspace{6mm}  \frac{1}{q_2} = \frac{1}{p_2} - \beta, \eqno(6.5)
$$

$$
K_{I,\alpha,\beta}(p_1, p_2) \stackrel{def}{=}
\sup_{ 0 \ne f \in L(p_1, p_2) } \left[\frac{|I^{\alpha, \beta}[f]|_{q_1, q_2}}{|f|_{p_1,p_2}} \right]. \eqno(6.6)
$$

\vspace{3mm}

{\bf Proposition 6.1.}  We state under necessary conditions (6.5)

$$
K_{I,\alpha,\beta}(p_1, p_2) = K_{I,\alpha}(p_1) \cdot  K_{I,\beta}( p_2). \eqno(6.7)
$$

\vspace{3mm}

{\bf Proof} is very simple. The {\it upper} estimate

$$
K_{I,\alpha,\beta}(p_1, p_2) \le K_{I,\alpha}(p_1) \cdot  K_{I,\beta}( p_2)
$$
may be proved analogously ones in the preprint  \cite{Ostrovsky18}, see also \ \cite{Ostrovsky10}.
 The {\it lower} bound for the variable $ K_{I,\alpha,\beta}(p_1, p_2) $ may be deduced by means
of choice of the factorable function

$$
f_0(x,y) = g_1(x) \cdot g_2(y),
$$
where the functions $ g_1(x),  \ g_2(y) $ are extremal functions for the correspondent one - dimensional problems. \par

\vspace{3mm}

{\bf Remark 6.1.} This circumstance, i.e. a phenomenon of factorability of the function $ K_{I,\alpha,\beta}(p_1, p_2) $ may
seem very strange, as long as in general case the individual operators of fractional integration and differentiation
are non commutative and $ |f|_{p_1, p_2} \ne |f|_{p_2, p_1}. $ \par

\vspace{3mm}

 Analogously may be grounded the following result. \\

 \vspace{3mm}

Denote

 $$
r_{\alpha, \beta}(p_1,p_2):= \sup_{0 \ne f \in B^{(\alpha, \beta)}_{p_1,p_2}}
\left\{ \frac{|D^{\alpha,\beta}[f]|_{p_1,p_2}}{||f||B^{(\alpha,\beta)}_{p_1,p_2}} \right\}, \eqno(6.8)
$$
where the double Besov's norm  $ ||f||B^{(\alpha,\beta)}_{p_1,p_2} $ denotes  the Bochner's composition of two Besov's norms
for the function $  f(x,y)  $ of two variables $  (x,y): $

$$
||f||B^{(\alpha,\beta)}_{p_1,p_2}  \stackrel{def}{=} || \ || \ f(\cdot,y) \ ||B_{x, p_1}^{(\alpha)} \ ||B_{y,p_2}^{(\beta)}.
$$

\vspace{3mm}

{\bf Proposition 6.2.} There holds under at the same necessary conditions (6.5)

$$
r_{\alpha, \beta}(p_1,p_2) = r_{\alpha}(p_1) \cdot r_{\beta}(p_2). \eqno(6.9)
$$

 \vspace{4mm}

 \section{Concluding remarks.}

 \vspace{3mm}

{\bf A. Weight  Riesz and Riemann-Liouville potential estimate.} \par

\vspace{3mm}

 An (linear)  operator $  U_{\alpha, \beta, \gamma}[f](x)   $ of a form

$$
U_{\alpha, \beta, \gamma}[f](x) = \frac{x^{-\gamma}}{\Gamma(\alpha)}  \ \int_0^x \frac{y^{-\beta} \ f(y) \ dy}{(x-y)^{\alpha - 1}}, \
x,y  \in (0,b), \ b = \const \le \infty \eqno(7.1)
$$
is called  weight  Riesz and Riemann-Liouville potential  operator, or generalized Cesaro-Hardy integral operator, or fractional integral. \par

 We refer here results about $  L_p \to L_q $ estimates for norm of this operator:

$$
|U_{\alpha, \beta, \gamma}[f]|_q \le V(\alpha, \beta, \gamma;p) \ |f|_p,
$$
when $ \alpha, \beta, \gamma  \in (0,1), \ \alpha + \beta + \gamma < 2,  \ \beta^2 + \gamma^2 > 0, $

$$
\frac{1}{q} =  \frac{1}{p} + (\alpha + \beta + \gamma - 2), \ 1 < p,q < \infty.
$$

and as ordinary

$$
V(\alpha, \beta, \gamma;p) \stackrel{def}{=} \sup_{0 \ne f \in L_p} \left[ \frac{|U_{\alpha, \beta, \gamma}[f]|_q }{|f|_p} \right], \ q = q(p).
$$

 This case is very different from the unweighed case, see, e.g.  \cite{Leoni1},  \cite{Lieb1}, \cite{Lieb2}, \cite{Lieb3}, \cite{Harboure1},
 \cite{Ostrovsky18}, \cite{Ostrovsky10} and so one.\par

 It is obtained in the aforementioned articles and books the $ L_p \to L_q $ sharp estimates for these operators, also in the
multidimensional case,  exact Grand Lebesgue norm  estimates etc. \par

 In detail, let us denote  $  \kappa = 2 - \alpha - \beta - \gamma,  $

$$
p_- = \frac{1}{1 - \beta}, \ p_+ = \frac{1}{2 - \alpha - \beta},
$$
and correspondingly

$$
q_- = \frac{1}{\alpha + \gamma - 1}, \ q_+ = \frac{1}{\gamma}.
$$

 Statement: if $ p \in (p_-, p_+],  $ or equally $ q \in [q_-, q_+),  $ then

$$
\frac{C_1(\alpha,\beta,\gamma)}{[ p - p_-]^{\kappa}}  \le V(\alpha, \beta, \gamma;p)    \le
\frac{C_2(\alpha,\beta,\gamma)}{[ p - p_-]^{\kappa}},  \ 0 < C_1 \le C_2 < \infty, \eqno(7.2)
$$
and $ V(\alpha, \beta, \gamma;p) = \infty  $ in other case. \par

\vspace{3mm}

{\bf B. Possible generalizations on metric measure spaces.} \par

\vspace{3mm}

 It is interest by our opinion to obtain our estimations, especially to derive the lower bounds, on the arbitrary  metric measure spaces
in the spirit of articles \cite{Maly1}, \cite{Nuutinen1} etc.;  where was applied the important notion of Riesz capacity. \par

 By definition, the Riesz potential of order $  \theta $ of a measurable function $ f:  X \to R, $  where the set $  X  $
is equipped by a distance function $  d = d(x,y) $ and by a Borelian non-trivial measure $  \tau,  $ is

$$
R_{(\theta)} [f](x) = \int_X \frac{f(y) \ \tau(dy)}{[\tau(B(x,y))]^{\theta}}, \eqno(7.3)
$$
see a recent work \cite{Fuglede1} and an article \cite{Gatto1}. \par


\vspace{4mm}

\begin{thebibliography}{99}

\vspace{3mm}

\bibitem{Adams1}
{\sc D.R.Adams, L.I. Hedberg.} {\it Function Spaces and Potential Theory.} Springer Verlag, Berlin, Heidelberg,
New York, 1996.

 \bibitem{Adams2}
{\sc D.R.Adams, R.J.Bagby.} {\it  Translation-dilation invariant estimates  for Riesz potential.} Indiana Univ. Math.
Journal. 1974, V.23, N 1, 1051 - 1067.

\bibitem{Adams3}
{\sc D.R. Adams. } {\it Choquet integrals in potential theory. } Publ. Mat. 42 (1998), 3-66.

\bibitem{Adams4}
{\sc D.R. Adams. } {\it On the existence of capacitary strong type estimates in } $  R^n,$
Arkiv f\"or Matematik, 14 (1976), 125-140.

\bibitem{Bapna1}
{\sc I. B. Bapna and Nisha Mathur.} {\it Application of Fractional Calculus in Statistics.}
Int. J. Contemp. Math. Sciences, Vol. 7, 2012, no. 18, 849-856

\bibitem{Bennet1}
{\sc Bennett C. and Sharpley R.} {\it Interpolation of operators.}  Orlando, Academic Press Inc.,1988.

\bibitem{Borla1}
{\sc Andrea Borla and Costen Protopoescu.} {\it  Nonparametric Estimation of the Fractional
Derivative of a Function Distribution. }
Internet publication, PDF, (2014).

\bibitem{Enikeeva1}
{\sc Farida Enikeeva. } {\it Adaptive minimax estimation of a fractional derivative.}
Statistics Probability Letters, {\bf 76},  (2006),  1441-1448.

\bibitem{Fiorenza1}
{\sc A.Fiorenza.} {\it Duality and reflexivity in grand Lebesgue spaces.}
Collectanea Mathematica (electronic version), {\bf 51}, 2, (2000), 131-148.

\bibitem{Fiorenza2}
{\sc A. Fiorenza and G.E. Karadzhov.} {\it Grand and small Lebesgue spaces and
their analogs.} Consiglio Nationale Delle Ricerche, Instituto per le
Applicazioni del Calcoto Mauro Picone, Sezione di Napoli, Rapporto tecnico n. 272/03, (2005).

\bibitem{Fuglede1}
{\sc Fuglede B. } {\it On the theory of potentials in locally compact spaces.} Acta. Math.
103,  (1960), 139-215.

\bibitem{Gatto1}
{\sc A.E. Gatto, C. Segovia and S. V'agi.} {\it On fractional differentiation and integration on spaces of homogeneous type.}
 Rev. Mat. Iberoamericana, {\bf 12}, (1996), 111-145.

\bibitem{Lieb3}
{\sc R.L.Frank and E.H.Lieb.} {\it Inversion Positivity and the sharp Hardy-Littlewood-Sobolev Inequality.}
Electronic Publications, arXiv:0904.4275v1 [math.FA] 27 Apr 2009.

\bibitem{Golubev1}
{\sc Golubev, G.K., Enikeeva, F.} (2001.) {\it On the minimax estimation problem of a fractional derivative.}
 Theory Probab. Appl. 46, 619-635.

\bibitem{Hardy0}
{\sc Hardy G.H.} {\it On some properties of integrals of fractional order. } Messenger. Math. 1917, V. 47 N 10, 145-150.

\bibitem{Hardy1}
{\sc Hardy G.H., Littlewood J.E.}
{\it  Some properties of fractional integrals.}
Proc. London Math. Soc., Ser. 2, (1928), V.24, 77-141.

\bibitem{Hardy2}
{\sc Hardy G.H., Littlewood J.E.}
{\it  Some properties of fractional integrals.} I.
Math. Zeitschrift,  (1928), V.27, N 4,  565-606.

\bibitem{Hardy3}
{\sc Hardy G.H., Littlewood J.E.}
{\it  Some properties of fractional integrals.} II.
Math. Zeitschrift,  (1932), V.34, N 34,  403-439.

\bibitem{Harboure1}
{\sc Harboure E., Macias R.A., Segovia C. }
{\it Boundedness of fractional operators on $ L(p) $ spaces with different weight. }
Trans. Amer. Soc., 1984, V.285 N 2, 629-647.

\bibitem{Hunt1}
{\sc R.A.Hunt. } {\it Developments Related to the A. E. Convergence of Fourier Series.} MAA Studies in Harmonic
Analysis. V.13, J. Math. Ass., (1998), p. 20-37.

\bibitem{Iwaniec1}
{\sc T.Iwaniec and C. Sbordone.} {\it On the integrability of the Jacobian under minimal hypotheses.}
Arch. Rat.Mech. Anal., 119, (1992), 129–143.

\bibitem{Iwaniec2}
{\sc T.Iwaniec, P. Koskela and J. Onninen.} {\it Mapping of finite distortion:
Monotonicity and Continuity}. Invent. Math. 144 (2001), 507-531.

\bibitem{Kozachenko1}
 {\sc Kozachenko Yu. V., Ostrovsky E.I.} (1985). {\it The Banach Spaces of
      random Variables of subgaussian type. }  Theory of Probab. and Math.
      Stat., (in Russian). Kiev, KSU, {\bf 32}, 43-57.

\bibitem{Leoni1}
{\sc G.Leoni.} {\it A first Course in Sobolev Spaces. Graduate Studies in Mathematics.} v. 105, AMS, Provi-
dence, Rhode Island, (2009).

\bibitem{Lieb1}
{\sc Lieb E.H, Loss M.} {\it Analysis.} Providence, Rhode Island, 1997.

\bibitem{Lieb2}
{\sc E.H.Lieb.} {\it Sharp constants in the Hardy-Littlewood-Sobolev and related inequalities.} Ann. of Math.,
(2), 118 (1983), no 2, 349-374.

\bibitem{Liouville1}
{\sc Liouville J.} {\it M\'emoire sur il' integration des \'equations  differ\'entielles a indices  fractionaires. }
J. l\'Ecole Roy. Polyt\'echn., 1835, V.15 Sect. 24, 17-54.

\bibitem{Maly1}
{\sc J. Maly and L. Pick. } {\it The sharp Riesz potential estimates in metric spaces,}
Indiana Univ. Math. J., 51 (2002), 251–268.

\bibitem{Nuutinen1}
{\sc Juno Nuutinen and Pilar Silvestre. } {\it The Riesz capacity in metric spaces.}
arXiv:1501.05746v1 [math.FA] 23 Jan 2015

\bibitem{Ostrovsky8}
{\sc E. Liflyand, E. Ostrovsky and L. Sirota.} {\it Structural properties of Bilateral Grand Lebesque Spaces.}
Turk. Journal of Math., {\bf 34,} (2010), 207-219.
TUBITAK, doi:10.3906/mat-0812-8

\bibitem{Ostrovsky1}
{\sc Ostrovsky E.I.} {\it  Exponential estimates for the random fields and its applications. } (1999), Moskow-Obninsk,
OINPE,  (in Russian).

\bibitem{Ostrovsky2}
{\sc Buldygin V.V., Mishtary D.Ch., Ostrovsky E.I., Puchalskii A.W.} {\it New Trends in  Probability Theory and Statistics. }
(1992),  VSP (Utrecht, Tokyo, New York).

\bibitem{Ostrovsky7}
{\sc E. Ostrovsky and L.Sirota.} {\it Moment Banach spaces: theory and applications.}
HAIT Journal of Science and Engineering, {\bf C}, Volume 4, Issues 1-2,
pp. 233 - 262, (2007).

\bibitem{Ostrovsky9}
{\sc E. Ostrovsky and L.Sirota.} {\it Well Posedness of the Problem of Estimation
Fractional Derivative for a Distribution Function.}
arXiv:1412.6829v1 [math.ST] 21 Dec 2014

\bibitem{Ostrovsky10}
{\sc E. Ostrovsky and L.Sirota.} {\it  Cesaro-Hardy operators on bilateral Grand Lebesgue Spaces.  }
arXiv:1307.5481v1 [math.FA] 20 Jul 2013

\bibitem{Ostrovsky11}
{\sc E. Ostrovsky and L.Sirota.} {\it  Riesz's and Bessel's operators in  bilateral Grand Lebesgue Spaces.  }
arXiv:0907.3321v1 [math.FA] 19 Jul 2009

\bibitem{Ostrovsky12}
{\sc E. Ostrovsky and L.Sirota.} {\it  Hardy-Littlewood inequalities for Riesz's potential.
Low bounds estimations for different powers.}
arXiv:0909.5663v1 [math.FA] 30 Sep 2009

\bibitem{Ostrovsky13}
{\sc E. Ostrovsky and L.Sirota.} {\it Weight Hardy-Littlewood inequalities for different powers.  }
arXiv:0910.5880v1 [math.FA] 30 Oct 2009

\bibitem{Ostrovsky16}
{\sc E. Ostrovsky and L.Sirota.} {\it Moment Banach spaces: theory and applications.} HAIT Journal of
Science and Engeneering, C, Volume 4, Issues 1-2, pp. 233-262, (2007).

\bibitem{Ostrovsky17}
{\sc E. Ostrovsky and L.Sirota.} {\it Boundedness of operators in bilateral Grand Lebesgue Spaces, with
exact and weakly exact constant calculation.}
 arXiv:1104.2963v1 [math.FA] 15 Apr 2011

\bibitem{Ostrovsky18}
{\sc E. Ostrovsky and L.Sirota.} {\it Multiple weight Riesz and Fourier transforms in bilateral anosotropic
Grand Lebesgue Spaces. } arXiv:1208.2392v1 [math.FA] 12 Aug 2012

\bibitem{Meerschaert1}
{\sc Mark Meerschaert, Jeff Mortensen, and Hans-Peter Scheffler.  } {\it Vector Gr\"unvald formula for fractional derivatives. }
 Internet electronic publication, 2014.

\bibitem{Miller1}
{\sc K. Miller and B. Ross.}  (1993) {\it An Introduction to Fractional Calculus and Fractional Differential Equations.}
 Wiley, New York.

\bibitem{Riemann1}
{\sc Riemann G.} {\it Versuch einer algemainen Auffasung der Integration und Differentiation.}
Gesamelte Math. Werke, Leipzig, Teubner Verlag, 1922-1923, 114-126.

\bibitem{Samko1}
{\sc S. G. Samko, A. A. Kilbas and O. I. Marichev.} {\it Fractional Integrals and Derivatives:
Theory and Applications.} Gordon and Breach Science Publishers, Yverdon, 1993.

\bibitem{Stein1}
{\sc E.M.Stein.} {\it Singular Integrals and Differentiability Properties of Functions.} Princeton University
Press, Princeton, (1992).

\bibitem{Stein2}
{\sc E.M.Stein and J.O. Str\"omberg.} {\it Behavior of maximal functions in $ R^n $ for large $ n.$ } Arkiv f\"or matematik.
Volume 21, no 1, May (1983); Published by Institut Mittag-Lefler; Djursholm, Sweden, p. 259-269.

\bibitem{Stein3}
{\sc E.M.Stein. } {\it The development of square finctions in the works of A.Zygmund.} Bull. Amer. Math. Soc.,
7, (1982), p. 359-376.

\bibitem{Stein4}
{\sc E.M.Stein.} {\it Maximal functions: Spherical means. } Proc. Nat. Acad, Sci. USA, 73 (1976), p. 2174-2196.

\end{thebibliography}
\end{document}